%% file: simplest_8tet-c2.tex
\documentclass[11pt]{amsart}

\usepackage{amsfonts, amstext, amsmath, amsthm, amscd, amssymb}
\usepackage{longtable, tabularx, array}
\usepackage{verbatim}
\usepackage{fullpage}

\newcommand{\done}{\hfill $\blacksquare$}

\renewcommand{\setminus}{{\smallsetminus}}

\newcommand{\deetwo}{\ensuremath{\mathbb{D}_2}}
\newcommand{\cyclic}{\ensuremath{\mathbb{Z}_2}}
\newcommand{\deefour}{\ensuremath{\mathbb{D}_4}}

\theoremstyle{plain}
\newtheorem{theorem}{Theorem}[section]

\newtheorem{lemma}[theorem]{Lemma}

\theoremstyle{definition}

\newtheorem{remark}[theorem]{Remark}

\newtheorem*{namedtheorem}{\theoremname}
\newcommand{\theoremname}{testing}

\title[The 500 simplest hyperbolic knots]{The 500 simplest hyperbolic knots}

\author[A.\ Champanerkar]{Abhijit Champanerkar}
\author[I.\ Kofman]{Ilya Kofman}
\author[T.\ Mullen]{Timothy Mullen}

\thanks{August 25, 2014}

\begin{document}

\begin{abstract}
We identify all hyperbolic knots whose complements are in the census
of orientable one-cusped hyperbolic manifolds with eight ideal
tetrahedra.  We also compute their Jones polynomials.
\end{abstract}

\maketitle

\section{Introduction}
Historically, complexity for knots has been measured by their minimal
crossing number; indeed, the tabulation of all 250 knots up to 10
crossings ushered in modern knot theory.  For a hyperbolic knot $K$, a
natural measure of its geometric complexity is the minimal number of
ideal tetrahedra needed to triangulate $S^3\setminus K$.  Hyperbolic
knots with geometric complexity up to six tetrahedra were found by
Callahan-Dean-Weeks \cite{CDW}, and extended to seven tetrahedra by
Champanerkar-Kofman-Paterson \cite{CKP}.  In this paper, we extend this
tabulation to eight tetrahedra, which completes the tabulation of knot
complements in the current census of orientable one-cusped hyperbolic
manifolds.

\ 

\begin{center}
\begin{tabular}{|c|cccccccc|c|}
\hline
Tetrahedra&1&2&3&4&5&6&7&8&$\leq$8\\
\hline
Manifolds&0&2&9&52&223&913&3388&12241&16828\\
\hline
Knots&0&1&2&4&22&43&129&301&502\\
\hline
\end{tabular}
\end{center}

\ 

Following the method outlined in \cite{CDW} and \cite{CKP}, we identified
all possible knot complements in the census of hyperbolic manifolds
provided in SnapPy \cite{snappy}, which is the version of SnapPea for
python developed by Culler and Dunfield.  Six Dehn fillings with trivial
first homology resisted efforts by SnapPy and Testisom \cite{testisom} to
simplify their fundamental groups, but all six were identified as
Seifert fibered using Matveev's 3-manifold Recognizer \cite{recognizer}.

Finding the corresponding knots is the challenge inherent in this
tabulation.  After exhaustive computer searches (similar to the ones
discussed in \cite{CKP}) using Knotscape \cite{knotscape} and various
generalized twisted torus knots, we were left with 51 unidentified
knot complements.  Because the geometric complexity is often only
mildly changed by small Dehn surgery, we continued our automated
search by randomly adding several full twists on subsets of strands of
various braid descriptions of low-complexity knots.  This left us with
24 knots to find ``by hand.''

SnapPy has a library of links (MorwenLinks) provided by Morwen
Thistlethwaite, and SnapPy identified complements of MorwenLinks after
drilling out various short geodesics from our knot complements.  We
then obtained link diagrams using Knotilus \cite{knotilus} or LinKnot
\cite{linknot}.  After obtaining surgery descriptions, we obtained the
knot diagrams using Kirby calculus either by hand or using the Kirby
Calculator \cite{kirbycalculator}.  Whenever possible, we simplified
these diagrams using Knotscape, Slavik Jablan simplified others for us
using his software, and finally we used Culler and Dunfield's new
python module Spherogram \cite{spherogram} to simplify the knot
diagrams.

\begin{remark} \rm
As discussed in \cite{BK}, many of the simplest hyperbolic knots are
Lorenz, even though few Lorenz knots have low crossing number.  The
following table gives updated numbers of known Lorenz knots among the
simplest hyperbolic knots.

\ 

\begin{center}
\begin{tabular}{|c|cccccccc|c|}
\hline
Tetrahedra&1&2&3&4&5&6&7&8&$\leq$8\\
\hline
Knots&0&1&2&4&22&43&129&301&502\\
\hline
Lorenz knots&0&0&1&2&12&23&70&141&249\\
\hline
\end{tabular}
\end{center}
\end{remark}

\ 

To identify all possible knot complements in the census, 
we used the following universal upper bounds on the Dehn filling
coefficients that may result in $S^3$.
We now show that there is a basis for homology of the
 cusp torus and a universal bound on $p$ and $q$ such that
whenever $p$ or $q$ is out of these bounds, the $(p,q)$--curve
on the maximal cusp has length greater than 6. Using the 6-Theorem
(\cite{agol2000, lackenby2000}), Dehn filling this curve results in a
3-manifold with an infinite, word-hyperbolic fundamental group, hence
not homeomorphic to $S^3$.  Therefore, to check for Dehn fillings
which result in $S^3$, one need only check a fixed set of Dehn filling
coefficients for all one-cusped manifolds.  Note that exceptional
surgeries can occur for arbitrarily long $(p,q)$--curves if there is no
restriction on the basis; e.g., there are knot complements with
exceptional Dehn fillings (lens spaces) using the standard basis with
arbitrarily long filling curves.

We take the shortest curve basis on the cusp torus, in which the
meridian is the shortest curve and the longitude is the second
shortest curve.  Note that SnapPy picks this basis by default for any
census manifold.  For a given manifold (e.g. a knot complement), one
can switch to such a basis using the SnapPy command
\verb|M.set_peripheral_curves(`shortest')|.

The idea can be summarized as follows. Taking the shortest curve basis on the
cusp torus, together with the lower bound on the length of the meridian on the
maximal cusp, imply that the parallelogram spanned by the meridian
and longtiude cannot be too thin, and hence has a lower bound on its
area. This gives an upper bound on how many such parallelograms
can be contained in a ball of radius 6, giving a bound on $p$ and $q$.
The more precise statement and proof follows from results in
\cite{cm2001,lm2013}.  Here we give the statement and its proof for completeness.

\begin{lemma}
\label{pqbounds}
For a complete orientable one-cusped hyperbolic 3-manifold $M$ with
cusp torus $T$, there is a basis $\mathcal{M}, \mathcal{L}$ of
$H_1(T)$ such that in the maximal cusp, the length of the curve
$p\mathcal{M} + q \mathcal{L}$ is greater than 6 if $|p| > 7 $ or $|q|
> 3$.
\end{lemma}

{\bf Proof:} Let the basis $\mathcal{M}, \mathcal{L}$ of $H_1(T)$ be
the shortest curve basis on a (fixed) cusp torus.
Without loss of generality, we can assume that $\mathcal{M}, \mathcal{L}$ 
are represented by the complex numbers $m \in \mathbb{R}^+$ and 
$\ell = a + ib,\ b > 0$ respectively. At the maximal cusp  
 $ m \geq 1$ (see e.g. \cite{adams2002}). Hence $|\ell| \geq m \geq 1$. 

Since $m$ and $\ell$ are shortest curves on the cusp torus, $|m \pm
\ell| \geq |\ell|$. Using the cosine law, this implies
$\displaystyle{|a|=|\ell\; \cos \theta| \leq \frac{m}{2}}$,
where $\theta$ is the angle between $m$ and $\ell$, hence 
$$b^2 = |\ell \; \sin \theta|^2 \geq |\ell|^2 - \frac{m^2}{4} \geq 
m^2 - \frac{m^2}{4} \geq \frac{3m^2}{4} 
\implies b \geq \frac{m\sqrt{3}}{2}. $$
Cao and Meyerhoff \cite{cm2001} show that the cusp area which equals $
b\cdot m \geq 3.35$, hence 
$$ \frac{2}{\sqrt{3}}\;b^2 \geq bm \geq 3.35 \implies b > 1.7.$$
Now $|pm + q\ell| = |pm + q(a + ib)| = |(pm+qa) + i(qb)| \leq 6$ 
implies that $|pm+qa| \leq 6$ and $ |qb| \leq 6$. Since $b > 1.7$, 
$|q| \leq 3$. Since $-6 \leq pm+qa \leq 6$, $|qa| \leq 3m/2$ and $m >1$ 
$$ -6 -\frac{3m}{2} \leq pm \leq 6 + \frac{3m}{2} \implies 
-6 -\frac{3}{2} \leq p \leq 6 + \frac{3}{2} 
\implies |p| \leq 7. $$
\done

\newpage

\subsection*{Acknowledgments}
As indicated above, this project relied heavily on various computer
tools developed over many years by others.  We are especially grateful
to Marc Culler and Nathan Dunfield for discussions and help with
aspects of SnapPy, which has become an essential tool to study
hyperbolic manifolds.  We also thank Culler and Dunfield for help with
Spherogram to simplify knot diagrams.  We also thank Nathan Dunfield
for help with Lemma \ref{pqbounds}.  We thank Slavik Jablan for
simplifying assorted knot diagrams.  We thank Saul Schleimer for
noticing that we inadvertently omitted two knot complements in an
earlier version of this paper.  We thank the anonymous referee for
providing Table 6 to us.  We are grateful for support by an NSF STEM
grant, Simons Foundation collaboration grants, and PSC-CUNY grants.

\section{Tables}

We follow \cite{CDW} for our notation and convention. Table 1 gives
a list of knots whose complements can be decomposed into 8 ideal
hyperbolic tetrahedra. We use the notation $k8_m$ to indicate the
$m^{th}$ knot in the list of knots made from $8$ tetrahedra. These
knots are sorted in increasing volume. When there are multiple knots
with the same volume, we then sort by decreasing length of systole, the
shortest closed geodesic in the complement. There are two pairs of
knots, ($k8_{53}$, $k8_{54}$ and $k8_{300}$, $k8_{301}$), which have
the same volume and the same length of systole.  In this case we have
ordered them by decreasing length of the second shortest closed
geodesic in the complement (computed using SnapPy command 
\verb|M.dual_curves()[1].complete_length().real|).

In Table 1 we also provide a description of a knot, either as a
Rolfsen census knot \cite{Rolfsen}, Hoste-Thistlethwaite census knot
\cite{HT}, or as a generalized twisted torus knot as described below.
The 59 knots which could not be described in one of the above ways are
denoted ``See Below,'' 
and we describe these knots in Table 2 using the Dowker-Thistlethwaite code.
%

In Tables 3 and 4, we update the descriptions for the $7$--tetrahedral census
knots.  In the published version \cite{CKP}, the last number was
deleted in the DT code for many knots.  Also, $k7_{80}$ was
incorrectly described in Table 1 of \cite{CKP}; it is the Perko knot,
$10_{161}=10_{162}$.  These errors were corrected in the ArXiv version
\cite{CKParxiv}.  Below, we provide new generalized twisted torus knot
descriptions which were not known before, and corrected or simplified
Dowker-Thistlethwaite codes for the other $7$--tetrahedral census
knots.

For all but one knot, we were able to compute the Jones polynomials of
the census knots, which are given in Table 5.  We also computed Jones
polynomials for two of the three knots not computed in \cite{CKP},
which are attached at the end of Table 5.

The anonymous referee kindly provided Table 6 with Conway symbols for census
knots which are not twisted torus knots.  As the referee pointed out,
``From the Conway symbols the structure of knots is clearly visible
(rational, Montesinos, etc.), and this info cannot be recognized from
the other notations... For some of the knots from this list I have not
succeeded in obtaining a `standard' Conway symbol (without single $-1$
crossings), so these knots are denoted by an asterisk. Some of them
can probably be written in `standard' form.''  



\subsection*{Generalized twisted torus knots}
Let $\delta_n = \sigma_1\cdots\sigma_{n-1}$ and $\delta^{-1}_n =
\sigma^{-1}_1\cdots\sigma^{-1}_{n-1}$ denote elements in the braid
group $B_n$.  Note these are not inverse braids; instead, they are mirror
images. 
The generalized twisted torus knot $T(r_1,s_1,\ldots,r_k,s_k)$ denotes
the closure of the braid $\delta_{r_1}^{s_1}\,\delta_{r_2}^{s_2}\cdots
\delta_{r_k}^{s_k}$.  In the tables below, our convention is that
$r_1>\ldots >r_k$.  (Lorenz knots are exactly T--knots with
$r_1<\ldots <r_k$ and all $s_i>0$ \cite{BK}.)

See \cite{ttk} for a discussion of various generalizations of twisted
torus knots.  In particular, in \cite{CDW} and \cite{CKP} the twisted torus
knot $T(p,q,r,s)$ denotes $s$ full twists on $r$ strands, which in our
notation would be denoted by $T(p,q,r,rs)$.

The knots are determined up to mirror image.  Note that in our
notation the mirror image of $T(p,q,r,s)$ is denoted by
$T(p,-q,r,-s)$, as in \cite{CDW} and \cite{CKP}.  This differs from \cite{ttk}
when $s$ is not a multiple of $r$ because our $\delta^{-1}_n$ is the
braid inverse of $\bar\delta_n=\sigma_{n-1}\cdots\sigma_1$ \cite{ttk}.


\ 

Here is a brief description of the table columns:
\begin{description}
\item[\textbf{Census}] This gives the number of the manifold in Thistlethwaite's 
8-tetrahedral census included in SnapPy; e.g. $S^3-k8_1$ is $t00017$.
\item[\textbf{Volume}]This is the hyperbolic volume of knot complement.
\item[\textbf{C-S}]This is the Chern-Simons invariant of the knot complement. 
\item[\textbf{Sym}]This is the group of isometries of the knot
  complement. In the table $0$ denotes the trivial group, \cyclic\
  denotes the cyclic group of order 2, \deetwo\ denotes the dihedral
  group of order 4 (Klein four) and \deefour\ denotes the dihedral
  group of order 8.
\item[\textbf{SG}]This is the length of the systole (shortest closed geodesic) in the knot
  complement.
\item[\textbf{Description}] Here we give a possible description of the
  knot using either Rolfsen's table of knots with 10 or fewer
  crossings, or Hoste-Thistlethwaite's table of knots from 11 to 16
  crossings, or as a generalized twisted torus knot as above.
\item[\textbf{Degree}]The first integer gives the lowest degree and the second integer gives the highest degree of the Jones polynomial.
\item[\textbf{Jones polynomial}]An entry of the form $(n,m)$ and
  $a_0+a_1+ \ldots +a_{n-m}$ corresponds to the polynomial $a_0t^n + a_1t^{n+1}+ \ldots + a_{n-m}t^m$.
\item[\textbf{N}]When the knot is from Rolfsen's table or
  Hoste-Thistlethwaite's table, this is the minimal crossing number of
  the knot.  For other knots, this is the crossing number of its
  diagram, which may not be minimal. 
\item[\textbf{DT code}]This is the Dowker-Thistlethwaite code for the knot.
\end{description}

\newpage


\begin{longtable}{|l||l|l|l|l|l|l|}

\caption{Knot complements with 8 tetrahedra}\\
  \hline
    & Census & Volume & C-S & Sym & SG & Description \\
  \hline 
  \endfirsthead

 \caption[]{\emph{continued}}\\
 \hline
   & Census & Volume & C-S & Sym & SG & Description \\
 \hline 
 \endhead

 \hline
 \endlastfoot

 \input{tables/table-8tet-c.txt}
\end{longtable}

\newpage

 \begin{longtable}{|l||l|l|p{4in}|}
 \caption{Dowker-Thistlethwaite codes for 59 census knots}\\
  \hline
   & Census & N & DT code \\
  \hline 
  \endfirsthead

 \caption[]{\emph{continued}}\\
 \hline
   & Census & N & DT code \\
   \hline 
 \endhead

 \hline
 \endlastfoot

 \input{tables/table-dt-8tet-c2.txt}
\end{longtable}

\newpage

\begin{longtable}{|l||l|l|>{\raggedright\arraybackslash}p{4in}|}
 \caption{Generalized twisted torus 7--tet census knots}\\
  \hline
   & Census & Description \\
  \hline 
  \endfirsthead

 \caption[]{\emph{continued}}\\
 \hline
  & Census & Description \\
   \hline 
 \endhead

 \hline
 \endlastfoot

 \input{tables/table-ttk-7tet-c2.txt}

 \end{longtable}

\begin{longtable}{|l||l|l|>{\raggedright\arraybackslash}p{4in}|}
 \caption{Revised DT codes for 7--tet census knots}\\
  \hline
   & Census & N & DT code \\
  \hline 
  \endfirsthead

 \caption[]{\emph{continued}}\\
 \hline
  & Census & N & DT code \\
   \hline 
 \endhead

 \hline
 \endlastfoot

 \input{tables/table-dt-7tet-c2.txt}
\end{longtable}

\newpage

 \begin{longtable}{|l||l|l|p{4in}|}

 \caption{Jones polynomials of census knots }\\
  \hline
  & Census & Degree& Jones polynomial \\
  \hline 
  \endfirsthead

 \caption[]{\emph{continued}}\\
 \hline
  & Census & Degree & Jones polynomial \\
 \hline 
 \endhead

 \hline
 \endlastfoot

 \input{tables/table-jp-c2.txt}
\end{longtable}

\ 

 \begin{longtable}{|l||l|l|p{4in}|}

 \caption{Conway notation for non-twisted-torus census knots}\\
  \hline
  & Knot & Conway notation \\
  \hline 
  \endfirsthead

 \caption[]{\emph{continued}}\\
 \hline
  & Knot & Conway notation \\
 \hline 
 \endhead

 \hline
 \endlastfoot

 \input{tables/table-conway.txt}
\end{longtable}

\bibliography{references.bib}
\bibliographystyle{plain}

\end{document}

%% file: tables/table-8tet-c.txt
%
%
$k8_{1}$ & 17 & 3.60046726 & -0.20486947  & \deetwo & 0.041  & $14a12741 $ \\ \hline
$k8_{2}$ & 27 & 3.60953917 & 0.05102925  & \deetwo & 0.035  & $15a54894 $ \\ \hline
$k8_{3}$ & 110 & 3.64226693 & -0.07514350  & \cyclic & 0.014  & $T(5,21,2,2) $ \\ \hline
$k8_{4}$ & 146 & 3.64611011 & -0.14759254  & \cyclic & 0.011  & $T(5,22,3,3) $ \\ \hline
$k8_{5}$ & 324 & 3.94544935 & -0.07096598  & \cyclic & 0.072  & $T(3,19,2,4) $ \\ \hline
$k8_{6}$ & 423 & 4.04447600 & -0.06469886  & \cyclic & 0.029  & $T(7,23,3,2) $ \\ \hline
$k8_{7}$ & 434 & 4.04647351 & -0.22135822  & \cyclic & 0.025  & $T(7,25,2,2) $ \\ \hline
$k8_{8}$ & 550 & 4.37490171 & 0.14531989  & \cyclic & 0.026  & $T(6,13,2,2) $ \\

%% file: tables/table-dt-8tet-c2.txt
$k8_{9}$ & t565 & 28 & -22, 42, -44, 46, -48, 50, -52, 54, -56, 2, -40, 18, -16, 14, -12, 10, -8, 6, -4, -20, -38, 36, -34, 32, -30, 28, -26, 24 \\ \hline
$k8_{10}$ & t577 & 29 & 18, -20, 22, -24, 26, -28, 50, -48, 46, -44, 42, -40, 38, -36, 34, 52, 30, 10, -8, 6, -4, 2, -58, 56, -54, 32, -12, 14, -16 \\ \hline
$k8_{20}$ & t1039 & 21 & -30, 6, 18, 24, -36, 22, -34, 20, -32, 4, 14, 10, -38, -40, -42, -2, -16, -12, -8, -26, -28 \\ \hline
$k8_{22}$ & t1142 & 28 & -20, 22, -38, 24, -36, 48, 32, -52, 54, -42, 40, 4, 8, 34, -10, -50, 28, 26, -6, -2, 56, -16, 14, -30, 12, -46, 44, 18 \\ \hline
$k8_{24}$ & t1235 & 29 & -28, 50, -16, 38, -22, 24, -58, 40, 52, -8, 36, 54, -56, 14, 32, 46, 12, 26, 20, 6, -2, 48, 30, 44, 42, 4, 18, -10, -34 \\ \hline
$k8_{28}$ & t1342 & 34 & -38, 54, 16, -14, 60, -28, -44, -56, 40, 2, -50, 34, 64, -62, -10, 26, 46, -66, 20, 52, -4, 8, 58, -24, 36, 68, -18, 42, 6, -12, 30, -32, 48, 22 \\ \hline
$k8_{39}$ & t1779 & 23 & -28, -14, 40, -10, -22, 42, -6, 38, -26, -44, -8, -20, -32, 36, -2, -16, 18, -46, 30, 4, 12, 24, -34 \\ \hline
$k8_{44}$ & t1901 & 29 & 10, -18, -30, 34, 44, -50, 4, 40, 52, 12, -58, -8, 36, 38, 54, -48, -56, -42, 26, 16, -28, 24, -22, -32, 20, 2, -14, -6, 46 \\

%% file: tables/table-ttk-7tet-c2.txt
$k7_{10}$ & v321 & T(9,8,7,-10) \\ \hline 
$k7_{24}$ & v595 & T(9,7,6,-10) \\ \hline 
$k7_{44}$ & v1126 & T(7,6,5,-9) \\ \hline 
$k7_{51}$ & v1392 & T(7,8,5,2) \\ \hline 
$k7_{55}$ & v1547 & T(7,9,5,2) \\ \hline 
$k7_{59}$ & v1690 & T(9,7,5,4) \\ \hline 
$k7_{67}$ & v1839 & T(9,10,8,4) \\ \hline 
$k7_{109}$ & v2900 & T(7,12,6,-4) \\ \hline 
$k7_{110}$ & v2925 & T(4,6,3,3,2,3) \\ \hline 
$k7_{115}$ & v3105 & T(4,7,3,5) \\ \hline

%% file: tables/table-dt-7tet-c2.txt
$k7_{11}$ & v329 & 23 & -30, -28, -24, -20, 38, 22, 26, 44, -34, -36, 10, -6, 12, -4, -2, -46, -42, -18, 40, 8, 16, -32, 14 \\ \hline
$k7_{26}$ & v656 & 22 & 12, 14, -38, 36, -26, 32, 2, 4, -6, -34, 28, -44, 20, -10, 42, -22, 8, -18, 16, 24, 30, 40 \\ \hline
$k7_{41}$ & v1063 & 18 & 14, -10, 18, -34, -4, -20, -22, 28, 6, 2, -12, 26, 32, -36, 16, 24, -8, -30 \\ \hline
$k7_{54}$ & v1434 & 18 & -10, -36, -24, -26, -2, 28, -30, 4, 8, -12, 14, -32, -34, -18, 20, -22, -6, 16 \\ \hline
$k7_{65}$ & v1811 & 17 & -12, -14, -16, -18, -28, -24, -2, -4, -6, 22, 30, 34, -10, 32, -8, 20, 26 \\ \hline
$k7_{74}$ & v1971 & 21 & -4, -32, 12, -14, 16, -24, 26, -40, -28, -36, -34, 30, -6, -42, 38, -2, -22, -20, -18, -10, 8 \\ \hline
$k7_{77}$ & v2001 & 17 & 10, 20, -12, 2, -28, -14, -6, 26, -24, 22, 4, -30, 32, -34, -8, 18, -16 \\ \hline
$k7_{84}$ & v2257 & 22 & -34, 36, -14, -30, 42, -6, 20, 38, -4, 12, 18, 2, -44, -8, -40, -26, 24, -22, 16, 10, -28, 32 \\

%% file: tables/table-jp-c2.txt
$k8_{1}$ & 17 & (-2,12) & $1-1+2-2+2-2+2-2+2-2+2-2+1-1+1$ \\ \hline
$k8_{2}$ & 27 & (1,16) & $1-1+2-2+2-2+2-2+2-2+2-2+2-1+1-1$ \\ \hline
$k8_{3}$ & 110 & (41,78) & $1+0+1+0+1+0+0+0+0+0+0+0+0+0+0+0+0+0+0+0+0+0+0-1+0-1+0+0+0+0+0+0+0+0+0+0-1+1$ \\ \hline
$k8_{4}$ & 146 & (45,85) & $1+0+1+0+1+0+0+0+0+0+0+0+0+0+0+0+0+0+0+0+0+0+0+0-1+0-1+0+0+0+0+0+0+0+0+0+0+0+0+1-1$ \\ \hline
$k8_{5}$ & 324 & (20,43) & $1+0+1+0+0+0+0+0+0+0+0+0+0+0+0+0+0+0+0+0+0-1+1-1$ \\ \hline
$k8_{6}$ & 423 & (68,122) & $1+0+1+0+1+0+1+0+0+0+0+0+0+0+0+0+0+0+0+0+0+0+0+0+0-1+0-1+0-1+0+0+0+0+0+0+0+0+0+0+0+0-1+1-1+1+0+0+0+0+0+0+0+1-1$ \\ \hline
$k8_{7}$ & 434 & (73,131) & $1+0+1+0+1+0+1+0+0+0+0+0+0+0+0+0+0+0+0+0+0+0+0+0+0+0-1+0-1+0-1+0+0+0+0+0+0+0+0+0+0+0+0+0+0+1-1+1-1+0+0+0+0+0+0+0+0-1+1$ \\ \hline
$k8_{8}$ & 550 & (31,57) & $1+0+1+0+1+0+1-1+1-1+1-1+1-1+1-2+1-2+1-2+2-2+2-2+1-1+1$ \\

%% file: tables/table-conway.txt
$k8_1$ & $14a12741$ & $12\,2$ \\  \hline
$k8_2$ & $15a54894$ & $13\,2$ \\  \hline
$k8_{54}$ & $14n6006$ & $9,2\,1,-2$ \\  \hline
$k8_{64}$ & $12a722$ & $9\,1\,2$ \\  \hline
$k8_{65}$ & $13a4874$ & $10\,3$ \\  \hline
$k8_{70}^*$ & $16n246032$ & $122^*-1:-1:-1:-1.2\,0.-2\,1\,0.-2\,0$ \\ \hline
$k8_{115}$ & $11a343$ & $7\,1\,3$ \\  \hline
$k8_{116}$ & $12a1166$ & $8\,4$ \\

%% file: simplest_8tet-c2.bbl
\begin{thebibliography}{10}

\bibitem{adams2002}
Colin~C. Adams.
\newblock Waist size for cusps in hyperbolic 3-manifolds.
\newblock {\em Topology}, 41(2):257--270, 2002.

\bibitem{agol2000}
Ian Agol.
\newblock Bounds on exceptional {D}ehn filling.
\newblock {\em Geom. Topol.}, 4:431--449, 2000.

\bibitem{BK}
Joan Birman and Ilya Kofman.
\newblock A new twist on {L}orenz links.
\newblock {\em J. Topology}, 2(2):227--248, 2009.

\bibitem{CDW}
Patrick~J. Callahan, John~C. Dean, and Jeffrey~R. Weeks.
\newblock The simplest hyperbolic knots.
\newblock {\em J. Knot Theory Ramifications}, 8(3):279--297, 1999.

\bibitem{cm2001}
Chun Cao and G.~Robert Meyerhoff.
\newblock The orientable cusped hyperbolic {$3$}-manifolds of minimum volume.
\newblock {\em Invent. Math.}, 146(3):451--478, 2001.

\bibitem{ttk}
Abhijit Champanerkar, David Futer, Ilya Kofman, Walter Neumann, and Jessica~S.
  Purcell.
\newblock Volume bounds for generalized twisted torus links.
\newblock {\em Math. Res. Lett.}, 18(6):1097--1120, 2011.

\bibitem{CKP}
Abhijit Champanerkar, Ilya Kofman, and Eric Patterson.
\newblock The next simplest hyperbolic knots.
\newblock {\em J. Knot Theory Ramifications}, 13(7):965--987, 2004.

\bibitem{CKParxiv}
Abhijit Champanerkar, Ilya Kofman, and Eric Patterson.
\newblock The next simplest hyperbolic knots, arXiv:math/0311380.

\bibitem{spherogram}
Marc Culler and Nathan~M. Dunfield.
\newblock {Spherogram, a computer program to simplify planar diagrams}.
\newblock {\tt http://\allowbreak snappy.\allowbreak computop.\allowbreak
  org/\allowbreak spherogram.\allowbreak html}.

\bibitem{snappy}
Marc Culler, Nathan~M. Dunfield, and Jeffrey~R. Weeks.
\newblock {SnapPy, a computer program for studying the geometry and topology of
  3-manifolds}.
\newblock {\tt http://\allowbreak snappy.\allowbreak computop.\allowbreak org}.

\bibitem{knotilus}
Ortho Flint and Stuart Rankin.
\newblock {Knotilus}.
\newblock {\tt http://knotilus.math.uwo.ca}.

\bibitem{testisom}
Derek Holt and Sarah Rees.
\newblock {TESTISOM}.
\newblock {Available by anonymous ftp at {\tt ftp.ncl.ac.uk:\allowbreak
  /pub/local/nser}}.

\bibitem{knotscape}
Jim Hoste and Morwen Thistlethwaite.
\newblock {Knotscape 1.01}.
\newblock {\tt http://www.math.utk.edu/\allowbreak
  $\sim$morwen/knotscape.html}.

\bibitem{HT}
Jim Hoste, Morwen Thistlethwaite, and Jeff Weeks.
\newblock The first 1,701,936 knots.
\newblock {\em Math. Intelligencer}, 20(4):33--48, 1998.

\bibitem{linknot}
Slavik Jablan and Radmila Sazdanovic.
\newblock {LinKnot}.
\newblock {\tt http://www.mi.sanu.ac.rs/\allowbreak vismath/linknot}.

\bibitem{lackenby2000}
Marc Lackenby.
\newblock Word hyperbolic {D}ehn surgery.
\newblock {\em Invent. Math.}, 140(2):243--282, 2000.

\bibitem{lm2013}
Marc Lackenby and Robert Meyerhoff.
\newblock The maximal number of exceptional {D}ehn surgeries.
\newblock {\em Invent. Math.}, 191(2):341--382, 2013.

\bibitem{recognizer}
Sergei Matveev.
\newblock {3--Manifold Recognizer}.
\newblock {\tt http://www.matlas.math.csu.ru}.

\bibitem{Rolfsen}
Dale Rolfsen.
\newblock {\em Knots and links}.
\newblock Publish or Perish Inc., Berkeley, Calif., 1976.
\newblock Mathematics Lecture Series, No. 7.

\bibitem{kirbycalculator}
Frank Swenton.
\newblock {Kirby Calculator}.
\newblock {\tt http://community.middlebury.edu/\allowbreak
  $\sim$mathanimations/kirbycalculator}.

\end{thebibliography}
